\documentclass[11pt]{article}

\usepackage{amssymb,amsmath,amsthm,amsfonts,fullpage}
\usepackage[dvips]{graphicx,color}

\usepackage[mathscr]{eucal}

\newtheorem{theorem}{Theorem}
\newtheorem{lemma}{Lemma}

\def\D{\mathbb{D}}
\def\R{\mathbb{R}}
\def\P{\mathbb{P}}

\begin{document}
\title{Approximation by power series with $\pm 1$ coefficients}
\date{\today}
\author{C. S{\.\i}nan G\"unt\"urk} 
\maketitle

\begin{abstract}
In this paper we construct certain type of 
near-optimal approximations of a class of 
analytic functions in the unit disc by power series with two distinct
coefficients. More precisely,
we show that if all the coefficients of the power series
$f(z)$ are real and lie in $[-\mu,\mu]$ where $\mu < 1$, 
then there exists a power series $Q(z)$ with coefficients in $\{-1,+1\}$
such that $|f(z)-Q(z)| \to 0$  at the rate 
$e^{-C/|1-z|}$ as $z \to 1$ non-tangentially inside the unit disc. 
A result by Borwein-Erd\'elyi-K\'os shows that 
this type of decay rate is best possible. The special
case $f \equiv 0$ yields a near-optimal 
solution to the ``fair duel'' problem of Konyagin.
\end{abstract}

\section{Introduction}

This paper is motivated by the ``fair duel'' problem which the author
heard from S.~Konyagin \cite{Konyagin}. The problem is the following:
There are two duellists
$X$ and $Y$ who will shoot at each other (only one at a time) using
a given $\pm 1$ sequence $q = (q_n)_{n \geq 0}$ which
specifies whose turn it is to shoot at time $n$. The shots are
independent and identically distributed random variables with outcomes
hit or miss. Each shot 
hits (and therefore kills) its target with small unknown 
probability $\epsilon$, which is arbitrary but fixed throughout the duel. The
``fair duel'' problem is to find an
ordering $q$, which is independent of $\epsilon$, and 
is as fair as possible in the sense that the probability
of survival for each duellist is as close to $1/2$ as possible. We measure
the fairness of an ordering $q$ by its {\em bias} function $B_q(\epsilon)$, 
defined to be
$$
B_q(\epsilon) :=  \P\{X \mbox{ survives}\} - 
\P\{Y \mbox{ survives}\},
$$
and ask that $B_q(\epsilon) \to 0$ as $\epsilon \to 0$ as fast as 
possible.\footnote{The problem makes sense only when we ask $q$ to be 
universal, i.e., independent of $\epsilon$. Otherwise, for any 
$\epsilon \leq 1/2$, the bias can be made zero by mapping the fractional 
$\beta$-expansion of the number 
$1/(2\epsilon)$ in the basis $\beta = 1/(1-\epsilon)$ to an ordering
\cite{Konyagin}.}

\par
It is elemantary to calculate the bias in terms of $q$.
Given that the duel is not
over before time $n$, which happens with probability $(1-\epsilon)^n$,
the probability that $Y$ is shot at time $n$ is equal to
$\epsilon$ if it is the turn of $X$ to shoot, and $0$ otherwise.
By symmetry of the problem, we find
$$ \P\{Y \mbox{ dies at time } n\} - 
\P\{X \mbox{ dies at time } n\} = \epsilon q_n (1-\epsilon)^n,
$$
where
we have assumed that $q_n=+1$ labels the turn of $X$ and $q_n=-1$ labels the
turn of $Y$. Summing over n, we obtain
$$ 
B_q(\epsilon) = \epsilon \sum_{n=0}^\infty q_n (1-\epsilon)^n.
$$

\par
At first, it may appear as the best ordering should be to simply alternate
between $X$ and $Y$, i.e., to set $q_n = (-1)^n$, for which
$B_q(\epsilon) = \epsilon/(2-\epsilon) = \Theta(\epsilon)$. However, 
this naive option is quickly
ruled out as for instance the $4$-periodic sequence given by
$q_0 = 1$, $q_1 = -1$, $q_2 = -1$, $q_3 = 1$ yields 
$B_q(\epsilon) = \epsilon^2/(1+(1-\epsilon)^2) = \Theta(\epsilon^2)$. 
Continuing in this fashion, it is tempting to 
think that the Thue-Morse sequence on the
alphabet $\{-1,+1\}$  (see, e.g. \cite{automatic})
might perhaps be the optimal sequence. For the Thue-Morse sequence, one has
$$B_{\mathrm{TM}}(\epsilon) = \epsilon \prod_{n=0}^\infty 
\left (1-(1-\epsilon)^{2^n} \right),
$$
where the infinite product $ \prod
\left (1-z^{2^n} \right) = \sum q_n z^n$ can in fact be taken as the definition
of this sequence. 
It is not difficult to show that 
there is a positive constant $c > 0$ such that 
$B_{\mathrm{TM}}(\epsilon) = \Omega(e^{-c(\log \epsilon)^2})$.
(See Section \ref{decay_TM} for a short derivation.)

\par
It turns out that one can do much better. 
One special outcome of this paper will be the construction of universal 
orderings $q$ for which $B_q(\epsilon)
= O(e^{-c/\epsilon})$ where $c>0$ is an absolute constant. 
In fact, we shall prove the following more general result:

\begin{theorem}\label{main_thm}
Let $0 \leq \mu < 1 \leq M < \infty$ be arbitrary and
$\mathscr{R}_M := \{z \in {\mathbb C} : |1-z| \leq M (1-|z|) \}$.
There exist constants $C_1:=C_1(\mu, M) > 0$  and $C_2:=C_2(\mu, M) > 0$
such that for any power series 
$$f(z) = \sum_{n=0}^\infty a_n z^n, 
~~~~~ a_n \in[-\mu,\mu], ~~ \forall n,$$
there exists a power series with $\pm 1$ coefficients, i.e., 
$$Q(z) = \sum_{n=0}^\infty q_n z^n,
~~~~~ q_n \in \{-1,+1\}, ~~ \forall n,$$ 
which satisfies
\begin{equation}
 |f(z)  - Q(z)| < C_1 e^{-C_2/|1-z|}
\end{equation}
for all $z \in \mathscr{R}_M \setminus \{1\}$.
\end{theorem}

\par
\begin{figure}[t]
\begin{center}
\includegraphics[height = 3in]{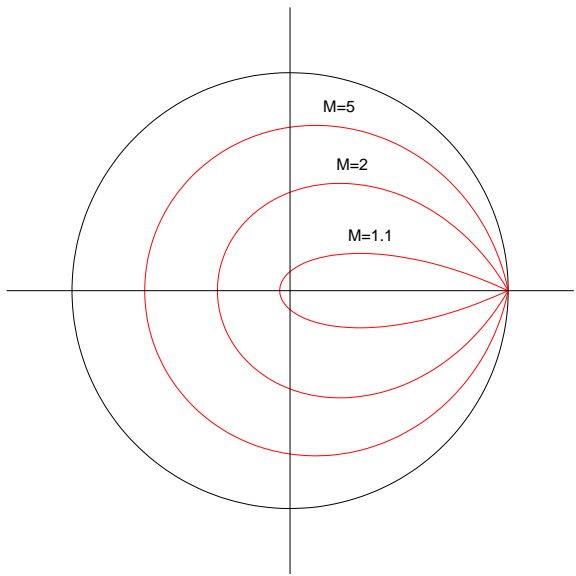}
\end{center}
\caption{\label{fig1} 
The boundary of $\mathscr{R}_M = 
\{z \in {\mathbb C} : |1-z| \leq M (1-|z|) \}$
for three different values of $M$.}
\end{figure}

Figure \ref{fig1} depicts the boundary of the 
set $\mathscr{R}_M$
for $M=1.1$, $M=2$ and $M=5$ along with the unit circle, which
is the limit as $M \to \infty$. Note that $\mathscr{R}_1 = [0,1]$. 
An alternative description
of this set can be given in the polar coordinates by 
$r \leq \exp(-\cosh^{-1}(1+(1-\cos\theta)/(M^2-1)))$ where we have assumed
that $\cosh^{-1}$ is given its positive value.

\par
In this generalized framework, the fair duel problem only corresponds to
approximating the zero function $f \equiv 0$ by power series with 
$\pm 1$ coefficients. For this case, we can set $\mu = 0$ and $M=1$.

\par
Note that in general it would be unrealistic to expect close
approximations by power series with $\pm 1$ coefficients at arbitrary
points inside the unit disc.
One extreme example is the point $z=0$. At this point, $|f(0)|=|a_0| 
\leq \mu < 1$ whereas $|Q(0)|=1$. In addition, it is not possible to have
a bound of the type 
$e^{-C/|a-z|}$ near any point $|a|<1$ as this
would imply $f\equiv Q$.

\par
It is clear that the point $z=1$ can be replaced by $z=-1$ by
considering  $\tilde a_n = (-1)^n a_n$ as input and setting 
$q_n = (-1)^n \tilde q_n$.
Also, the theorem extends to the case of arbitrary complex coefficients
$a_n$ such that $|\Re(a_n)| \leq \mu$ and $|\Im(a_n)| \leq \mu$ if
we allow $q_n = \pm 1 \pm i$.

\par
Apart from constants, the result of Theorem \ref{main_thm} is optimal 
by the following theorem of Borwein-Erd\'elyi-K\'os.

\begin{theorem}[Borwein-Erd\'elyi-K\'os {\cite[Thm 5.1]{Borwein}}]
\label{borwein_lower_bound}
There are absolute constants $c_1>0$ and $c_2 > 0$ such that
for any analytic function $F$ defined on the open unit disc $\D$ that satisfies
$$|F(z)| \leq \frac{1}{1 - |z|}, \;\;\;\;\;\;\;\; z \in \D,$$
one has, for any $\alpha \in (0,1]$,
$$ |F(0)|^{c_1/\alpha} \leq e^{c_2/\alpha} \sup_{x \in [1-\alpha,1]}
|F(x)|.$$
\end{theorem}

\par
To see the optimality, it suffices to set $F(z) = (f(z)-Q(z))/2$. Then
for any choice of $a_n$ and $q_n$ in their given range, $F$ satisfies the
hypothesis of Theorem \ref{borwein_lower_bound} and 
moreover $1 > |F(0)| \geq (1-\mu)/2$.

\section{Proof of Theorem \ref{main_thm}}

Given any sequence $(a_n)_{n \geq 0}$ 
that takes values in the interval $[-\mu,\mu]$, 
we will construct a $\pm 1$ sequence $(q_n)_{n \geq 0}$ 
via the following algorithm. Let $h = (h_k)_{k \geq 1}$,
be a sequence of real numbers (to be specified later) which satisfies
\begin{equation} \label{h_norm}
\mu+ \sum_{k=1}^\infty |h_k| \leq 2. 
\end{equation}
We shall call such a sequence $\mu$-admissible.
For $n = 0,1,2,\dots$, let
\begin{eqnarray}
w_n & = & \sum_{k=1}^n h_k v_{n-k} \label{w_eq} + a_n,\\
q_n & = & \mbox{sign}(w_n), \\
v_n & = & w_n - q_n, 
\end{eqnarray}
with the convention $\mbox{sign}(0) = 1$. Note also that $w_0 = a_0$. 

An easy induction argument gives that $|v_n| \leq 1$ for all $n$: It is
true for $n=0$ since $v_0 = a_0 - \mbox{sign}(a_0)$.
Assume  that $|v_k| \leq 1$ for all $k < n$. Then (\ref{h_norm}) and
(\ref{w_eq}) yield $|w_n| \leq 2$ and 
therefore $v_n = w_n - \mbox{sign}(w_n) \in [-1,1]$.

Since
$$ a_n - q_n = v_n - \sum_{k=1}^n h_k v_{n-k}, $$ 
we have the relation
\begin{equation}
\sum_{n=0}^\infty a_n z^n - \sum_{n=0}^\infty q_n z^n
= \left (1 - \sum_{k=1}^\infty h_k z^k \right)
\left (\sum_{n=0}^\infty v_n z^n \right)
\end{equation}
for $|z| < 1$.
Hence, using the boundedness of $v_n$, we obtain the estimate
\begin{equation}
\big |f(z) - Q(z) \big | 
\leq   \left |1 - \sum_{k=1}^\infty h_k z^k \right |  (1-|z|)^{-1},
\label{bound1}
\end{equation}
and the problem is reduced to finding a $\mu$-admissible sequence $h$
such that the power series 
$$H(z) = 1 - \sum_{k=1}^\infty h_k z^k$$ 
decays very rapidly to $0$ as $z \to 1$. We shall pick a particular
sequence $h$ as follows: Let $\sigma$ be a positive
integer. Let 
$$c_\sigma := \frac{\sinh (\pi{\sigma^{-1/2}})}{\pi{\sigma^{-1/2}}},$$
and let $h_k := h^{(\sigma)}_k$ be defined via
\begin{equation}
H(z) := H_\sigma(z) :=
1 - c_\sigma z - 2 c_\sigma \sum_{n=1}^\infty 
\frac{(-1)^n}{\sigma n^2 + 1} \, z^{\sigma n^2 + 1}.
\end{equation}
Below we shall estimate the size of $H_\sigma(z)$ for $|z| < 1$.
But first we have to show that 
(\ref{h_norm}) is satisfied. 

\par
We say that
$\sigma$ is $\mu$-admissible if $h^{(\sigma)}$ is $\mu$-admissible.
For any $\mu < 1$, this is obviously the case for all large values of $\sigma$
since $c_\sigma \to 1$  as $\sigma \to  \infty$. 
It is easy to find the $\mu$-admissible values of $\sigma$ explicitly.
We have the formula 
(e.g., \cite[p. 268]{henrici})
\begin{equation}
 1 + 2\, \omega^2 \sum_{n=1}^\infty 
\frac{1}{n^2 + \omega^2} = \pi \omega \coth(\pi\omega),
\end{equation}
which implies that 
\begin{eqnarray}
\sum_{k=1}^\infty |h^{(\sigma)}_k|
& = & c_\sigma\left(1 + 2 \sigma^{-1} \sum_{n=1}^\infty 
 \frac{1}{n^2 + \sigma^{-1}} \right )  \nonumber \\
& = & c_\sigma \pi \sigma^{-1/2} \coth(\pi\sigma^{-1/2}) \nonumber \\
& = & \cosh(\pi\sigma^{-1/2}).
\end{eqnarray}
Hence $\sigma$ is $\mu$-admissible if and only if
$$ \sigma > \pi^2/
\log^2 \left (2-\mu + \sqrt{(1-\mu)(3-\mu)} \right ).$$
We note that the smallest attainable 
value of $\sigma$ is $6$ and is $\mu$-admissible 
for $\mu = 2-\cosh(\pi/\sqrt{6})$. It is easily seen that
$c_\sigma$ is a decreasing function of $\sigma$; therefore 
$c_\sigma \leq c_6 < 1.3$ for all $\sigma$.
On that other hand, as $\mu \to 1$, 
the lower bound for $\sigma$ behaves as $\pi^2/2(1-\mu)$. 

\par
Of course the function $H_\sigma(z)$ was not chosen arbitrarily. 
First, we claim 
\begin{equation}
\label{zero_val}
\lim_{z \to 1} H_\sigma(z) = H_\sigma(1) 
= 0.
\end{equation}
The first equality follows from the uniform convergence of $H_\sigma(z)$ on the
closed unit disc and the second equality
follows straight from the formula (e.g., \cite[p. 271]{henrici})
\begin{equation}
\sum_{n=1}^\infty \frac{\cos n\xi}{n^2+\omega^2}
= \frac{\pi}{2\omega}\, \frac{\cosh(\pi-\xi)\omega}{\sinh \pi\omega}
- \frac{1}{2 \omega^2}
\end{equation}
for $\xi = \pi$ and $\omega = \sigma^{-1/2}$.

\par
In order to estimate $H_\sigma(z)$ around $z=1$, we note that
for $|z| < 1$, 
\begin{equation}
\label{deriv_theta}
H'_\sigma(z) = - c_\sigma \left ( 1 + 2 \sum_{n=1}^\infty (-1)^n 
z^{\sigma n^2} \right)
= -c_\sigma \varTheta_4(0,z^\sigma),
\end{equation}
where 
$$\varTheta_4(\omega,z) := \sum_{n=-\infty}^{\infty} (-1)^n z^{n^2} 
e^{2in\omega}$$
is the fourth theta-function of Jacobi. Hence, (\ref{zero_val}) and
(\ref{deriv_theta}) now provide us with
\begin{eqnarray} \label{bias_theta}
H_\sigma(z) 
=  c_\sigma \int_z^1 \varTheta_4(0,s^\sigma) \,ds.
\end{eqnarray}
We take the path of integration to be the straight line segment
$[z,1]$ connecting $z$ to $1$ to obtain the estimate
\begin{equation}
|H_\sigma(z)| \leq c_\sigma 
|1-z| \sup_{s \in [z,1)} |\varTheta_4(0,s^\sigma)|.
\end{equation}

\par
It is not difficult to estimate 
$|\varTheta_4(0,s^\sigma)|$ near $s=1$. One way of doing this 
is the classical method of applying the
Poisson summation formula to the modulated Gaussian
$f_\lambda(u) := e^{-\pi \lambda u^2} e^{i \pi u}$ where $\lambda > 0$
and $u \in \R$. Since $\widehat{f_\lambda}(\xi) = 
e^{-\pi (\xi-\frac{1}{2})^2 / \lambda}$, we have the identity
\begin{equation}\label{theta} 
\varTheta_4(0,e^{-\pi \lambda}) =
\sum_{n = -\infty}^{\infty} (-1)^n e^{-\pi \lambda n^2}
=  \frac{1}{\sqrt{\lambda}} 
\sum_{n = -\infty}^{\infty} e^{-\pi  (n-\frac{1}{2})^2 /\lambda },
\end{equation}
which then extends
to any complex $\lambda$ with $\Re(\lambda) > 0$ by analytic continuation,
using the principle branch of $\sqrt{\lambda}$. 
Clearly the dominating terms for this last expression
are given by $n=0$ and $n=1$. We replace 
$\lambda$ by $\lambda \sigma / \pi$ in (\ref{theta}) and
set $s = e^{-\lambda}$ so that
\begin{eqnarray}
\left | \varTheta_4(0,s^\sigma) \right |
& \lesssim & 
\frac{1}{|\lambda \sigma|^{1/2}}
\sum_{n = -\infty}^{\infty} 
\Big | e^{-\frac{\pi^2}{\lambda \sigma }(n-\frac{1}{2})^2 }
\Big | \nonumber \\
& \lesssim & \frac{1}{|\lambda \sigma|^{1/2}}
\left [\sum_{n = 0}^{\infty} 
\Big | e^{-\frac{\pi^2}{\lambda \sigma }}\Big |^n \right ]
\Big | e^{-\frac{\pi^2}{4 \lambda \sigma }}\Big |
\nonumber \\
& \lesssim & 
\label{theta-bound}
\frac{1}{|\lambda \sigma|^{1/2}}
\left [ 1+\frac{\sigma}{\Re\!\left(\frac{1}{\lambda}\right)} \right ]
e^{-\frac{\pi^2}{4 \sigma}\Re\left(\frac{1}{\lambda}\right)},
\end{eqnarray}
where in the last step we have used the inequality 
$(1-e^{-x})^{-1} < 1 + x^{-1}$ which is valid for all $x > 0$.
(We note that by $A \lesssim B$ we mean $A \leq C B$ for
an absolute positive constant $C$.)
Of course, this upper bound works best if we choose $\lambda$
with $-\pi \leq \Im(\lambda) < \pi$ since
$\Re\!\left(\frac{1}{\lambda}\right) = \Re(\lambda)/|\lambda|^2$.

\par
Let us finish the proof of Theorem \ref{main_thm}. First we note that
we do not lose any generality if we exclude the values of 
$z$ in the compact set
$\mathscr{K}_M := \left ( \mathscr{R}_M \cap \{z : \Re(z) \leq 0\} \right )
\cup \{z : |z| \leq 1/2\}$ since for any choice of $a_n$ and $q_n$
in their given range, the function $f(z) - Q(z)$ is bounded on 
$\mathscr{K}_M$ by a constant that only depends on $M$ and on this set
Theorem \ref{main_thm} does not yield a
bound better than a constant anyway. Hence we assume
that $z \in \mathscr{R}_M \setminus \mathscr{K}_M$.

\par Next is a simple lemma that we will employ to finish our analysis.
\begin{lemma}\label{simple_lemma}
Let $z \in \mathscr{R}_M \setminus \mathscr{K}_M$.
There exists an absolute constant $c>0$ such that if
$s = e^{-\lambda} \in [z,1)$ where $\lambda$ is chosen such that
$-\frac{\pi}{2} \leq \Im(\lambda) \leq \frac{\pi}{2}$, 
then $\Re\!\left(\frac{1}{\lambda}\right) \geq \frac{c}{M|1-s|}$.
\end{lemma}

\begin{proof}[Proof of Lemma \ref{simple_lemma}]
Let us write $s = |s|e^{-i\phi}$ so that
$\phi := \Im(\lambda)$ and $\Re(\lambda) = \log \frac{1}{|s|}$.
The fact that $z \in \mathscr{R}_M \setminus \mathscr{K}_M$
and $s \in [z,1)$ imply that $|s| \geq c$ for some absolute constant $c > 0$. 
Hence we have 
$$\log \frac{1}{|s|} \lesssim  |s| \log \frac{1}{|s|} \leq 1-|s|
\leq |1 - s|.$$
On the other hand, we also have
$$|\phi| \lesssim |\sin \phi| \lesssim |s| |\sin \phi | 
= |\Im(s)| = |\Im(1-s)| \leq |1-s| \leq M(1-|s|) \leq M \log \frac{1}{|s|},$$
where the fourth inequality relies on the fact that $\mathscr{R}_M$ is a
convex set. This, together with the previous inequality, imply
$$ \frac{1}{\Re\!\left(\frac{1}{\lambda}\right)} = 
\log \frac{1}{|s|} + \frac{\phi^2}{\log \frac{1}{|s|}}
\lesssim M |1-s|.
$$
This proves the lemma.
\end{proof}

Now, using the lemma, we have the bound 
\begin{equation}
\left [ 1+\frac{\sigma}{\pi^2 \,\Re\!\left(\frac{1}{\lambda}\right)} \right ]
e^{-\frac{\pi^2}{4 \sigma}\Re\left(\frac{1}{\lambda}\right)}
\lesssim \sigma M e^{-\frac{2C}{\sigma M |1-s|}}.
\end{equation}
where $C$ is an absolute positive constant.
On the other hand, we have 
$$|\lambda| \geq \Re(\lambda) = 
- \log |s| \geq 1 - |s| \geq \frac{1}{M} |1-s|,$$ 
so that
$$\frac{1}{|\lambda|^{1/2}}  
\lesssim \sigma^{1/2} M 
\left( \frac{C}{\sigma M |1-s|}\right)^{1/2} 
\lesssim \sigma^{1/2} M e^{\frac{C}{\sigma M |1-s|}}.$$
Hence we obtain the desired estimate 
$$ |f(z)-Q(z)| \lesssim \sigma^{3/2} M^3 e^{\frac{-C}{\sigma M |1-z|}} $$
which concludes the proof of Theorem \ref{main_thm}.
\hfill $\Box$

\section{Remarks} \label{remarks}

\subsection{Explicit upper bounds for the optimal decay of bias} 
For the special case $f \equiv 0$ that corresponds to the
fair duel problem, we can set $\mu = 0$, $M=1$ 
and choose any $\sigma \geq 6$. Let us call the resulting sequence
$q^{(\sigma)}$. 
Moreover, we are only interested in approximation of $f(z)$ for real values
of $z$. It is easy to check from the proof of Theorem \ref{main_thm} 
that we now have
$$ |B_{q^{(6)}}(\epsilon)| \lesssim \sqrt{\epsilon}\; 
e^{-\frac{\pi^2}{24 \epsilon}}. 
$$

\par
Below is the first $50$ values of $q^{(6)}_n$ computed using our
algorithm. For the compactness of presentation, 
we list it as a $\{0,1\}$ sequence rather than a $\pm 1$ sequence.

$$
\begin{array}{rl}
{\tt q^{(6)}} :  &
{\tt 
100101011010101001011
010010101011010010110
10100101...
}
\end{array}
$$
Note that it would be necessary to employ special numerical methods to
compute the terms for arbitrarily large $n$ due to 
the possibility of the accumulation of rounding errors. 

\par
It is interesting that the beginning of the
sequence $q^{(8)}$ bears a remarkable resemblance with
the Thue-Morse sequence:

$$
\begin{array}{rl}
\tt q^{(8)} :  &
{\tt 
1001011001101001
1001011001101001
100101100110010110...} \\
\tt TM:  &
{\tt 
1001011001101001
0110100110010110
011010011001011010...}
\end{array}
$$

\subsection{Decay of bias for the Thue-Morse sequence}\label{decay_TM}

For the completeness of our discussion we present below a short derivation
of the decay of bias for the Thue-Morse sequence. 
Let $N_\epsilon$ be the unique integer such that $1 \leq 2^{N_\epsilon}
\epsilon < 2$. We have
$$ 1 \geq \prod_{n = N_\epsilon}^\infty 
\left (1-(1-\epsilon)^{2^n} \right)
\geq \prod_{n = N_\epsilon}^\infty 
\left (1-e^{-2^n \epsilon} \right) 
\geq \prod_{n = 0}^\infty 
\left (1-e^{-2^n} \right) 
\gtrsim 1;
$$
hence it suffices to estimate the product of the first $N_\epsilon
= \left(\log_2 \frac{1}{\epsilon}\right)(1+o(1))$
terms. For this it suffices to use the simple inequality
$2^n \epsilon \geq 1-(1-\epsilon)^{2^n} \geq \epsilon$. We now have
$$
2^{-\frac{1}{2}\left(\log_2 \frac{1}{\epsilon}\right)^2(1+o(1))}
= \prod_{n = 0}^{N_\epsilon-1} 2^n \epsilon
\geq \prod_{n = 0}^\infty 
\left (1-(1-\epsilon)^{2^n} \right)
\gtrsim \prod_{n = 0}^{N_\epsilon-1} \epsilon
= 2^{-\left(\log_2 \frac{1}{\epsilon}\right)^2(1+o(1))}.$$

\subsection{Extensions}

The proof of Theorem \ref{main_thm} employed the boundedness of the sequence
$(v_n)$. It is possible to relax this condition by allowing
for mild (e.g. sub-exponential) growth of $|v_n|$. 
This generalization is one possible direction to seek
better approximations. 

Our algorithm was inspired by sigma-delta quantization; in fact,
the particular scheme that we have employed corresponds to an
``infinite-order'' limit of a family of 
schemes developed in \cite{exp_decay}.

\section*{Acknowledgements}

The author would like to thank Sergei Konyagin for introducing him
the fair duel problem and the many valuable 
discussions. This work was supported in part 
by the National Science Foundation Grant DMS-0219072.

\vskip 1 cm

\noindent 
C. S{\.\i}nan G\"unt\"urk \\
Courant Institute of Mathematical Sciences\\
251 Mercer Street \\
New York, NY 10012. \\
{\tt gunturk@cims.nyu.edu}


\begin{thebibliography}{12}

\bibitem{Konyagin}
S.~Konyagin, 
\newblock personal communication, 2003.

\bibitem{automatic}
J-P.~Allouche and J.~Shallit,
\newblock ``Automatic Sequences: Theory, Applications, Generalizations,''
\newblock Cambridge University Press, 2003.

\bibitem{henrici}
P.~Henrici,
\newblock  ``Applied and Computational Complex Analysis,''
\newblock Wiley, 1991.

\bibitem{Borwein}
P.~Borwein, T.~Erd\'{e}lyi and G.~K\'{o}s,
\newblock ``Littlewood-type problems on $[0,1]$,''
\newblock {\em Proc. London. Math. Soc.} 79 (1999), 22--46.


\bibitem{exp_decay}
C.~S. G{\" u}nt{\" u}rk,
\newblock ``One-Bit Sigma-Delta Quantization with Exponential Accuracy,''
\newblock {\em Comm. Pure Appl. Math.},  vol. 56, pp. 1608--1630, no. 11,
2003.

\end{thebibliography}
\end{document}